\documentclass[12 pt]{amsart}

\usepackage{amsmath}
\usepackage{amssymb} 
\usepackage{latexsym}

\newtheorem{theorem}{Theorem}[section]
\newtheorem{corollary}[theorem]{Corollary}
\newtheorem{lemma}[theorem]{Lemma}

\begin{document}  
 
\title[Symmetry of Arthur parameters under Aubert involution]{Symmetry of Arthur parameters under Aubert involution}

\author{Dubravka Ban}

\maketitle

\begin{center}
{\it Department of Mathematics,
Southern Illinois University\\ Carbondale, IL 62901, USA}

E-mail: {\it dban@math.siu.edu}

\end{center}

\bigskip

{\Small
  \noindent ABSTRACT.
We consider a nontempered $A$-parameter $\psi$ of $SO(2n+1, F)$ of a certain type
and the base point representation $\pi$ in the $A$-packet of $\psi$.
 Let $\hat{\pi}$ be the Aubert involution of $\pi$.
We compute explicitly the Langlands data of $\hat{\pi}$
and the $A$-parameter
$\hat{\psi}$  of $\hat{\pi}$. We investigate whether $\psi$ and 
$\hat{\psi}$ are symmetric. Although symmetry holds for large classes of
parameters, it does not hold in general.
}  
 
\bigskip

\section{Introduction}

This paper deals with large classes of nontempered representations of
the odd orthogonal group $SO(2n+1,F)$  over a $p$-adic field $F$.
These representations arise from considerations of $A$-parameters 
of a certain type (\ref{A1}).
In accordance with Arthur's conjectures \cite{A1,A2},
attached to each $A$-parameter is a finite set of equivalence classes of
irreducible admissible representations, called an $A$-packet.
There is, however, a natural way to associate to each $A$-parameter
a particular representation; we call it a base point.
We study effects of the duality operator on $A$-parameters via base points.
For a nontempered representation $\pi$ with the $A$-parameter $\psi$,
we compute explicitly the Langlands data of the dual representation
and the corresponding $A$-parameter. The proof relies on recent fundamental
developments by Jiang-Soudry, Harris-Taylor and Henniart.
It provides an interesting illustration of the Langlands-Arthur
 formalism.
Recall that
$A$-parameters and $A$-packets emerged from Arthur's work on the question of 
how nontempered representations should fit into the trace formula. 
There are very few examples
beyond tempered parameters for larger groups, where Arthur's formalism has been 
confirmed.

The duality operator is a generalization of the Zelevinsky involution.
The Zelevinsky involution is an operator defined on the 
Grothendieck group of the category of all smooth finite length
representations of the general linear group $GL(n,F)$ \cite{Zele}.
This involution has many important properties. It relates a discrete series
representation to the corresponding Langlands quotient.
The Zelevinsky involution on $GL(n,F)$ preserves unitarity. 
Furthermore, its action on $A$-parameters can be precisely defined, as follows.
Let
\[
  \psi : W_F \times SL(2, \mathbb{C}) \times SL(2, \mathbb{C}) \to GL(n, \mathbb{C})
\]
be an $A$-parameter of $GL(n,F)$. Here, $W_F$ denotes the Weil group of $F$.
Let $\pi$ be the representation of $GL(n,F)$ associated to $\psi$.
Denote by $\hat{\pi}$ the Zelevinsky involution of $\pi$ and by $\hat{\psi}$
the $A$-parameter of $\hat{\pi}$. Then \cite{Zele, MoeW, Tad:unitary},
\begin{equation}\label{A0}
     \hat{\psi}(w, x,y) = {\psi}(w,y,x).
\end{equation}
In other words, the Zelevinsky involution acts on $A$-parameters by interchanging
two copies of $SL(2, \mathbb{C})$. We say $\psi$ and $\hat{\psi}$ are symmetric.

The Zelevinsky involution allows generalizations to a connected reductive 
quasi-split algebraic group $G$ defined over $F$. 
Bernstein \cite{Be}, Schneider and Stuhler \cite{Sch-stu},
and Aubert \cite{Aub} have defined
duality operators on
the category of all smooth finite length representations of $G$
and on its Grothendieck group.
The duality operator sends an irreducible representation to an irreducible
representation. Other questions, related to important properties of the 
Zelevinsky involution, are still open. It is expected that the duality 
operator preserves unitarity, which seems to be very difficult to prove.
Barbasch and Moy in \cite{Bar-Moy} proved the conjecture for representations
with nonzero Iwahori-fixed vectors, using the Kazhdan-Lustzig parametrization
of such representations. 
Even more interesting is the question of the action of the involution on $A$-packets.
Barbasch conjectured that the duality operator 
sends an $A$-packet to an $A$-packet. 
If Barbasch's conjecture holds, we may consider the $A$-parameter associated to 
an $A$-packet and the $A$-parameter associated to the packet obtained by applying
the duality operator on the original packet.
This raises the question of the action
of the involution on $A$-parameters. It is conjectured that,
as for general linear groups,  the involution 
acts on $A$-parameters of $G$  
by interchanging
two copies of $SL(2, \mathbb{C})$.
Although the conjecture was known previously,
a precise statement is due to Hiraga \cite{Hiraga}.
 In a joint work with Zhang \cite{BanZh},
we proved that,
for a generic discrete series representation ${\pi}$ of $SO(2n+1,F)$,
the $A$-parameters of ${\pi}$ and $\hat{\pi}$ are symmetric.
 This has further consequences; for example, this justifies 
a generic discrete series representation of a Levi subgroup of $SO(2n+1,F)$
and its involution have the same $R$-group, as conjectured by Arthur
(cf. \cite{BAnnSci, BAlgebra}).

In this paper, we consider certain nontempered  $A$-parameters.
Let $\rho$ be an irreducible unitary supercuspidal representation of $GL(n,F)$
and $\sigma$ an irreducible supercuspidal generic representation of $SO(2\ell+1,F)$.
Assume $\rho$ is self-dual (equivalent to its contragredient).
Denote by $S_n$ the $n$-dimensional irreducible complex representation of
$SL(2,\mathbb{C})$.
Let $\phi$ be the L-parameter of $\rho$ \cite{H-T, Henn}.
According to the work of
Jiang and Soudry \cite{J-S}, we can find
the L-parameter of $\sigma$. It is of the form
$
     \bigoplus_{i \in A} \phi_i \otimes S_1.
$
Let $\pi$ be the base point representation
associated to the A-parameter 
\begin{equation}\label{A1}
    {\psi}= \phi \otimes S_{k} \otimes S_{2}
\oplus
     \bigoplus_{i \in A} \phi_i \otimes S_1 \otimes S_1,
\end{equation}
$k \geq 1$. 
Denote by $\hat{\pi}$ the Aubert involution of $\pi$ and by $\hat{\psi}$
the $A$-parameter of $\hat{\pi}$.
In accordance with the conjectures explained above, one may expect
\begin{equation}\label{A2}
    \hat{\psi}= \phi \otimes S_{2} \otimes S_{k}
\oplus
     \bigoplus_{i \in A} \phi_i \otimes S_1 \otimes S_1.
\end{equation}
If (\ref{A2}) holds, we say the $A$-parameters of $\pi$
and  $\hat{\pi}$ are symmetric.
Symmetry is equivalent to (\ref{A0}), or, colloquially, the condition that
the Aubert involution acts on the $A$-parameter of $\pi$ by interchanging
two copies of $SL(2,\mathbb{C})$.
We prove  that symmetry of the parameters depends on 
both
the parity of $k$ and  the 
point of reducibility $\alpha$ of the induced representation 
$i_{G,M}(\nu^\alpha \rho \otimes \sigma)$.
For example, if $i_{G,M}(\nu^\frac{1}{2} \rho \otimes \sigma)$
is reducible, then we have the following:
 if $k$ is even,
\[
  \hat{\psi} =  \phi \otimes S_1 \otimes S_{k+1} \oplus \phi \otimes S_1 \otimes S_{k-1}
     \oplus   \bigoplus_{i \in A} \phi_i \otimes S_1 \otimes S_1;
\]
if $k$ is odd, 
\[
   \hat{\psi} = \phi \otimes S_{2} \otimes S_{k} \oplus 
      \bigoplus_{i \in A} \phi_i \otimes S_1 \otimes S_1
\]
(Theorem \ref{Th1}). 
In other words, if  $k$ is odd, the parameters of $\pi$ and $\hat{\pi}$ 
are symmetric; if $k$ is even, they are not.
Similar relations, with the parity of $k$ interchanged,
hold for $\alpha = 0$ (Theorem \ref{Th2}). The case 
$\alpha = 1$ is slightly different (Theorem \ref{Th3}).
Let us point out that this does not imply the $A$-packets, with the 
corresponding $A$-parameters, are 
not symmetric under the Aubert involution. Indeed,  
the Aubert involution does not always send a base point to a base point
(see Remark 6.1; for results on $U_{E/F}(4)$, see \cite{Konno}).
Therefore, we do not prove or disprove the $A$-packets are symmetric
under the Aubert involution. Our work concerns base points and $A$-parameters.
This knowledge, however, is essential for understanding behavior of $A$-packets.
In addition, the results are exact and do not depend on
conjectures. Namely,
the hypotheses about Plancherel measures and unitarity of involution,
assumed in \cite{BAnnSci, BAlgebra, BanZh},
are  not assumed in this work.

The base point associated to an $L$-parameter is determined based on the work
of Jiang and Soudry \cite{J-S}. They deal with groups $SO(2n+1,F)$ and in this 
paper we consider the same series of groups.
In view of
the recent work by 
Cogdell, Kim, Piatetski-Shapiro and Shahidi \cite{CKP-SS},
we expect our methods
can be applied to other series of classical $p$-adic groups.

We now give a short summary of the paper. In Section 2, we recall some basic
definitions and properties of $L$-parameters and $A$-parameters.
In Section 3, we prove some technical lemmas 
 on Jacquet modules of parabolically induced representations.
The lemmas are needed in the rest of the paper and
the proofs of the main results (Theorems  \ref{Th1}, \ref{Th2} and \ref{Th3})
rely on considerations of 
Jacquet modules.
In Sections 4, 5 and 6, we study symmetry of the $A$-parameters 
under the Aubert involution.
Although symmetry holds for large classes of
parameters,  it does not hold in general.
As mentioned earlier, symmetry of the parameters depends on 
the 
point of reducibility $\alpha$ of the induced representation 
$i_{G,M}(\nu^\alpha \rho \otimes \sigma)$.
Each section is devoted to one of the cases $\alpha= \frac{1}{2}$, 0
and 1.

{\em Acknowledgment.} The conjecture that the involution 
acts on $A$-parameters by interchanging
two copies of $SL(2, \mathbb{C})$ was introduced to me by
Anne-Marie Aubert and Peter Schneider in Luminy, 2002.
James Arthur explained to me the importance of the conjecture 
during the Clay Mathematics Institute Summer School 
at the Fields Institute, 2003.
This paper has benefited from discussions with Dan Barbasch,
David Goldberg, Colette M{\oe}glin, Gordan Savin
and Freydoon Shahidi. I thank them all.

\section{Preliminaries}

In this section, we recall some basic
definitions and properties of $L$-parameters and $A$-parameters
and do some preliminary computation on $A$-parameters.
Let $F$ be a nonarchimedean local field of characteristic zero and 
$G$  a reductive group over $F$.
Let $P$ be a standard parabolic subgroup of $G$ with the Levi decomposition $P=MU$.

\subsection{Parabolic induction and segments}
If $\sigma$ is
a smooth representation of $M$, we denote by
$i_{G,M}(\sigma )$ the representation parabolically induced from $\sigma$.
For a smooth representation $\pi$ of $G$, $r_{M,G}(\pi)$ is normalized
Jacquet module of $\pi$ with respect to $M$ \cite{BeZ, Cass}.
For admissible representations 
$\rho_i$ of $GL(k_i,F)$,  $i=1,2,$  define
\[
    \rho_1 \times \rho_2 = i_{G,M}(\rho_1 \otimes \rho_2),
\]
where 
$M \cong GL(k_1,F) \times GL(k_2,F)$ is a standard Levi subgroup of 
$G= GL(k_1+k_2,F)$.
If $\rho$ and $\sigma$ are
admissible representations of $GL(k,F)$ and $SO(2\ell+1,F)$,  define
\[
    \rho \rtimes \sigma = i_{G,M}(\rho \otimes \sigma),
\]
where 
 $M \cong GL(k,F) \times SO(2\ell+1,F)$ is a standard Levi subgroup of $G=SO(2(k+\ell)+1,F)$
\cite{Tad:unitary}.

Define $\nu =|det|$.
Let $\rho$ be an irreducible supercuspidal
representation of $GL(k,F)$ and $m \leq n$  integers. The set
$[\nu^m\rho,\nu^n \rho] = \{\nu^m \rho, \nu^{m+1} \rho, \dots, \nu^n \rho \}$ is called
a segment \cite{Zele}. The induced representation 
$
   \nu^n \rho \times \nu^{n-1} \rho \times \cdots \times \nu^m \rho 
$
 has a unique
irreducible subrepresentation, which we denote by $\delta [\nu^m \rho,\nu^n \rho]$.

For a representation $\sigma$, we denote by  $\tilde{\sigma}$ the contragredient 
of $\sigma$.
If $\rho$ is an irreducible supercuspidal representation of $GL(k,F)$,
then 
$
    \widetilde{(\nu^\alpha \rho)} = \nu^{-\alpha} \tilde{\rho}.
$
For a segment $\Sigma =[\rho,\nu^n \rho]$, we define
$
     \tilde{\Sigma} = [\nu^{-n} \tilde{\rho}, \tilde{\rho}].
$
Then 
$
    \delta(\tilde{\Sigma}) = \widetilde{\delta(\Sigma)}
$
\cite{BeZ}.

\subsection{Grothendieck group and Aubert involution}
Let $R(G)$ be the Grothendieck group of the category of all smooth finite 
length representations of $G$.
For a smooth finite length representation $\pi $ of $G$,  we define
 $s.s.(\pi ) \in R(G)$ to be the sum of the
irreducible components of $\pi$, each component taken with the multiplicity
corresponding to its multiplicity in $\pi$. 
Let $\leq$  denote
the natural partial order on $R(G)$.
For smooth finite length representations $\pi_1$  and $\pi_2$,
 we write $\pi_1 \leq \pi_2$ if
$s.s.(\pi_1) \leq s.s.(\pi_2)$ in the Grothendieck group.

The Aubert duality operator $D_G$ is defined on the Grothendieck group \cite{Aub}.
If $\pi$ is an irreducible admissible representation of $G$, we define
$
      \hat{\pi} =\pm D_G(\pi),
$
taking the sign + or - so that $\hat{\pi}$ is a positive element in the
Grothendieck group. We call $\hat{\pi}$ the Aubert involution of $\pi$.
It follows from \cite{Aub} that $\hat{\pi}$ is an irreducible representation.

\subsection{Langlands classification for $SO(2n+1,F)$}
Suppose $\rho_i$ is a discrete series representation of $GL(n_i,F)$, 
$i=1, \dots ,k$ and $\alpha_1 \leq \cdots \leq \alpha_k < 0$ are 
real numbers.
Let $\tau$ be a tempered representation of  $SO(2\ell+1,F)$. 
Then the induced representation 
$
   \nu^{\alpha_1} \rho_1 \times \cdots \times \nu^{\alpha_k} \rho_k \rtimes \tau
$
has a unique irreducible subrepresentation, which we call the Langlands
subrepresentation and denote by 
$
  L_s(\nu^{\alpha_1} \rho_1 , \dots , \nu^{\alpha_k} \rho_k , \tau ).
$
For any irreducible admissible representation $\pi$ of $SO(2n+1,F)$, there exist
Langlands data $\nu^{\alpha_1} \rho_1 , \dots , \nu^{\alpha_k} \rho_k , \tau$
as above, such that 
$
  \pi = L_s(\nu^{\alpha_1} \rho_1 , \dots , \nu^{\alpha_k} \rho_k , \tau ).
$

\subsection{Langlands parameters and base points}
Let
$W_F$ be the Weil group of $F$.
We take $W_F \times SL(2, \mathbb{C})$ as the Weil-Deligne group \cite{Tat, Knapp}.
A Langlands parameter, or $L$-parameter, of $SO(2n+1,F)$ is 
a 
homomorphism
\[
  \phi : W_F \times SL(2, \mathbb{C})  \to Sp(2n, \mathbb{C})
\]
such that $\phi(W_F)$  consists of semi-simple elements in $Sp(2n, \mathbb{C})$
and the restriction of $\phi$ to $SL(2, \mathbb{C})$ is algebraic
\cite{Bor, Langl, Knapp}. The parameter $\phi$ is tempered if the image $\phi(W_F)$
is bounded.
Two $L$-parameters are equivalent if they are conjugate in $Sp(2n, \mathbb{C})$.
According to the Local Langlands Conjecture, 
each parameter $\phi$ should parametrize a finite set of equivalence classes
of irreducible admissible representations of $SO(2n+1,F)$,
called the $L$-packet of $\phi$.
Langlands proved the correspondence for real and complex groups
\cite{Langl}.
Precise description of $L$-packets for real groups is due to
Shelstad \cite{Shel:indis}.
The Local Langlands Conjecture for $GL(n,F)$ was proved 
by Harris and Taylor \cite{H-T}, and Henniart \cite{Henn}.
Jiang and Soudry in \cite{J-S}  defined a bijection
\[
       \phi \, \longleftrightarrow \, \pi= 
L_s(\nu^{\alpha_1} \rho_1 , \dots , \nu^{\alpha_k} \rho_k , \tau ), \quad \tau \text{ generic}
\]
 between the set of equivalence classes of
$L$-parameters of $SO(2n+1,F)$ and the set of equivalence classes of 
irreducible admissible representations 
$
\pi= 
L_s(\nu^{\alpha_1} \rho_1 , \dots , \nu^{\alpha_k} \rho_k , \tau ),
$
with $\tau$ generic. This bijection preserves $L$- and $\epsilon$-factors.
The representation $\pi$ is a member of the $L$-packet
of $\phi$ and
plays an important role. For example,
 if $\phi$ is tempered, then $\pi = \tau$ is generic, 
which confirms the conjecture 
formulated by Shahidi in
\cite{Shah:proof}
on existence of a generic member in each tempered packet.
We call the representation $\pi$ associated to $\phi$ by \cite{J-S}
the base point representation in the $L$-packet of $\phi$.

\subsection{Arthur parameters}
An Arthur parameter, or $A$-parameter, of $SO(2n+1,F)$ is 
a 
homomorphism
\[
  \psi : W_F \times SL(2, \mathbb{C}) \times SL(2, \mathbb{C}) \to Sp(2n, \mathbb{C})
\]
such that $\psi(W_F)$ is bounded and included in the set 
of semi-simple elements of $Sp(2n, \mathbb{C})$
and the restriction of $\psi$ to the two copies of $SL(2, \mathbb{C})$ is algebraic
\cite{A1, A2, Kim-Sh}.
In accordance with Arthur's conjectures, attached to each $A$-parameter
$\psi$ is a finite set of equivalence classes of irreducible
admissible representations, called the $A$-packet of $\psi$.
To any $A$-parameter $\psi$, Arthur associates an $L$-parameter $\phi_\psi$ by
\[ \phi_{\psi}(w,x)=\psi(w,x, \left(\begin{array}{cc} 
                                      |w|^{1/2} &   \\
                                                 & |w|^{-1/2}
                                    \end{array}
                              \right) ).
\]
Contrary to $L$-packets, $A$-packets need not to be disjoint. 
A representation $\pi$ may occur in more than one packet.
An $A$-parameter $\psi$ is called the $A$-parameter of $\pi$ 
if $\phi_\psi$ is the $L$-parameter of $\pi$.
This definition is justified by noticing that 
$ \psi \mapsto \phi_\psi$ is injective \cite{A2}.
If $\psi $ is an A-parameter, we may decompose 
it into a direct sum 
\[
       \psi = \bigoplus_{i=1}^k (\phi_i \otimes S_{m_i} \otimes S_{n_i}),
\]
where $m_i, n_i \in {\Bbb Z}^+$, 
$\phi_i$ is a continuous homomorphism such that
$\phi_i (W_F)$ consists of semisimple matrices and $S_m$ is the $m$ dimensional
irreducible complex representation of $SL(2,\Bbb C)$. 
Note that
\[ 
 \phi(w)\otimes S_{n}(\left(\begin{array}{cc} 
                   |w|^{1/2} &   \\
                            & |w|^{-1/2}
              \end{array}      
        \right) )
   =\bigoplus_{j = -(n-1)/2}^{(n-1)/2}\; 
          \phi(w)|w|^j.
\]
Therefore, for $\psi = \phi \otimes S_{m} \otimes S_{n}$, we have
\begin{equation}
\label{A}
     \phi_\psi = \bigoplus_{j = -(n-1)/2}^{(n-1)/2} \; 
          |\cdot|^j \phi \otimes S_{m}.
\end{equation}

Let $\sigma$ be an irreducible supercuspidal generic representation 
of $SO(2\ell+1,F)$  and let $\rho$ be an irreducible
supercuspidal unitary representation of $GL(k,F)$.
Let
$
     \bigoplus_{i \in A} \phi_i \otimes S_1
$
be the L-parameter of $\sigma$ defined in \cite{J-S}. 
Let $\phi$ be the L-parameter of $\rho$.
We consider 
\[
   \psi = \phi \otimes S_{m} \otimes S_2 \oplus 
      \bigoplus_{i \in A} \phi_i \otimes S_1 \otimes S_1,
\]
where $m \in \Bbb{Z}^+$.
The corresponding L-parameter $\phi_{\psi}$
is equal to
\[
    \phi_{\psi} = 
     |\cdot|^{\frac{1}{2}} \phi \otimes S_{m}  \oplus |\cdot|^{-\frac{1}{2}} \phi \otimes S_{m}
     \oplus
\bigoplus_{i \in A} \phi_i \otimes S_1
\]
and the base point representation attached to this $L$-parameter is
\begin{equation}\label{B}
    \pi = L_s( 
   \delta [\nu^{-\frac{m}{2}} \rho ,\nu^{\frac{m}{2}-1}\rho], \sigma).
\end{equation}

\section{Jacquet modules}

In the proofs of the main results (Theorems  \ref{Th1}, \ref{Th2} and \ref{Th3}),
we rely on considerations of 
Jacquet modules of parabolically induced representations.
In this section, 
we prove some technical lemmas on Jacquet modules we need in the rest of the paper.

The following lemma follows directly from \cite{BCan}, Corollary 4.3.

\begin{lemma}\label{Can3} 
Let $\sigma$ be an irreducible supercuspidal representation of $M$ and
 $\pi$ an irreducible subquotient of $i_{G,M}(\sigma)$. 
Then $\pi$ is a subrepresentation of $i_{G,M}(\sigma)$ if and only if 
$ \sigma \leq r_{M,G}(\pi)$.

\end{lemma}

\begin{lemma}\label{lemma}
Let $\pi$ be an admissible representation of $G=SO(2 \ell+1,F)$.
Let
\[
    \rho_1 \otimes \cdots \otimes  \rho_i \otimes \rho_{i+1} \otimes \cdots 
    \otimes \rho_k \otimes \sigma
\]
be an irreducible supercuspidal representation of a standard Levi subgroup $M$ of $G$.
Assume 
$
   r_{M,G}(\pi) \geq  \rho_1 \otimes \cdots \otimes  \rho_i \otimes \rho_{i+1} \otimes \cdots 
    \otimes \rho_k \otimes \sigma.
$

\begin{enumerate}

     \item[{(i)}] If $\rho_i \times \rho_{i+1} \cong \rho_{i+1} \times \rho_{i}$, then
$
   r_{M,G}(\pi) \geq  \rho_1 \otimes \cdots \otimes  \rho_{i+1} \otimes \rho_{i} \otimes \cdots 
    \otimes \rho_k \otimes \sigma.
$

    \item[{(ii)}] If $\rho_k \rtimes \sigma$ is irreducible, then
$
   r_{M,G}(\pi) \geq  \rho_1 \otimes \cdots \otimes  \rho_i \otimes \rho_{i+1} \otimes \cdots 
    \otimes \rho_{k-1} \otimes \tilde{\rho}_k \otimes \sigma.
$

\end{enumerate}

\end{lemma}

{\it Proof.} We prove (i). For (ii), the proof is similar. According to Lemma \ref{Can3}, 
\[
   \pi   \hookrightarrow
 \rho_1 \times \cdots \times  \rho_{i} \times \rho_{i+1} \times \cdots 
    \times \rho_k \rtimes \sigma
  \cong
 \rho_1 \times \cdots \times  \rho_{i+1} \times \rho_{i} \times \cdots 
    \times \rho_k \rtimes \sigma
\]
Lemma \ref{Can3} tells us  
$
   r_{M,G}(\pi) \geq  \rho_1 \otimes \cdots \otimes  \rho_{i+1} \otimes \rho_{i} \otimes \cdots 
    \otimes \rho_k \otimes \sigma.
$
\hfill   $\square$

For an ordered partition $(n_1, \dots, n_q = n)$ of $n$, denote by 
$
Sh_{(n_1, \dots, n_q)}
$
the set of all shuffles of sets
$
\{1,\dots,n_1\},\, \{n_1+1,\dots,n_2\}, \dots , \{n_{q-1}+1, \dots, n_q \}
$
\cite{BCan}.
(Suppose that $S_1,S_2, \dots ,S_q$ are disjoint ordered sets. A shuffle of the
sets $S_1,S_2, \dots ,S_q$ is a permutation $p$ of the set 
$S= S_1 \cup S_2 \cup \cdots \cup S_q$ which preserves the order on each of the sets
$S_k$.)
For $k \leq l \leq n$, define a permutation $z_{(k,l)}$ with
\[
z_{(k,l)}(j)=
     \begin{cases}
       j, & j<k \\
       k+l-j, & k \leq j \leq l; \\ 
      j, & j>l.
     \end{cases}
\]
If $k>l$, we define $z_{(k,l)}=1$.
Set 
$
      1_q = \underbrace {1, \dots, 1}_{q-times}
$ and
$
     -1_q = \underbrace {-1, \dots, -1}_{q-times}.
$
 Let 
$
    M \cong GL(n,F)^k \times SO(2\ell+1, F)
$
be a standard Levi subgroup of $G = SO(2(nk+\ell)+1, F)$.
Let
$
    \pi_0 = \rho_1 \otimes \cdots \otimes \rho_k \otimes \sigma
$ 
be an irreducible supercuspidal representation of 
$M$.
A permutation of the set $\{ 1, \dots, k \}$ acts on $\pi_0$ by permuting 
$
    \rho_1, \dots , \rho_k.
$
In addition, if 
$
   \epsilon = (\epsilon_1, \dots, \epsilon_k) \in \{ \pm 1 \}^k,
$
then $\epsilon$ acts on $\pi_0$ as
\[
   \epsilon \pi_0 = \rho_1^{\epsilon_1} \otimes \cdots \otimes 
     \rho_k^{\epsilon_k} \otimes \sigma, \quad  \text{ where }
   \rho_i^{\epsilon_i} = \begin{cases}
                          \rho_i, &if \, \epsilon_i=1,\\
                          \tilde{\rho_i}, &if \, \epsilon_i=-1.
                         \end{cases}
\]
The following lemma follows from \cite{BCan}. It describes Jacquet modules of 
representations induced from intermediate Levi subgroups.
 
\begin{lemma} 
\label{Can}
 Let 
$
    M \cong GL(n,F)^k \times SO(2\ell+1, F)
$
be a standard Levi subgroup of $G = SO(2(nk+\ell)+1, F)$ and 
$
    N \cong GL(pn,F) \times SO(2((k-p)n+\ell)+1, F)
$
 an intermediate Levi subgroup, $M < N < G$.
Let $\pi_0$ be an irreducible supercuspidal representation of 
$M$. Suppose that $\pi$ is a subquotient of
$i_{N,M}(\pi_0)$. Then
\[
    s.s.( r_{M,G} \circ i_{G,N}(\pi)) =
        \sum_{q=0}^{p}  
      Sh_{(p-q, p, k)} z_{(p-q+1, p)} (1_{p-q}, -1_q, 1_{k-p})\,
      r_{M,N}(\pi).
\]
\end{lemma}

\begin{lemma}\label{Can2}
 Let 
$
    M \cong GL(n,F)^k \times SO(2\ell+1, F)
$
be a standard Levi subgroup of $G = SO(2(nk+\ell)+1, F)$ and
$
    N \cong GL(pn,F) \times GL((k-p)n,F) \times SO(2\ell+1, F)
$
 an intermediate Levi subgroup, $M < N < G$.
Let 
$
\pi_0 = \rho_1 \otimes \cdots \otimes \rho_k \otimes \sigma
$ 
be an irreducible supercuspidal representation of 
$M$. 
If $\pi$ is a subquotient of
$i_{N,M}(\pi_0)$, then
\[
\begin{aligned}
    s.s.&( r_{M,G} \circ i_{G,N}(\pi)) \\
     &=
        \sum_{q=0}^{p}      
      \sum_{r=0}^{k-p}
      Sh_{(p-q, p, k-r, k)} z_{(p-q+1, p)} z_{(k-r+1, k)}
      (1_{p-q}, -1_q, 1_{k-p-r}, -1_r) 
      r_{M,N}(\pi).
\end{aligned}
\]

\end{lemma}

{\it Proof.} Let 
$
    N_1 \cong GL(pn,F) \times SO(2((k-p)n+\ell)+1, F)
$
and 
$
   \pi_1 = i_{N_1, N} (\pi).
$
Then, by Lemma \ref{Can}, 
\[
   \begin{aligned}
    s.s.&(r_{M,G} \circ i_{G,N}(\pi)) = r_{M,G} \circ i_{G,N_1}(\pi_1) \\
      &=         \sum_{q=0}^{p}      
        Sh_{(p-q, p,k)} z_{(p-q+1, p)}
      (1_{p-q}, -1_q, 1_{k-p}) 
      r_{M,N_1}(\pi_1) \\
   &=        \sum_{q=0}^{p}      
        Sh_{(p-q, p, k)} z_{(p-q+1, p)}
      (1_{p-q}, -1_q, 1_{k-p}) 
      \sum_{r=0}^{k-p}
      Sh^{(p)}_{(k-r, k)} z_{(k-r+1, k)}(1_{k-r}, -1_r) r_{M,N_1}(\pi_1),
  \end{aligned}     
\]
where $Sh^{(p)}_{(k-r, k)}$ denotes the subset of $Sh_{(k-r, k)}$ consisting of
permutations which keep $\{1, \dots, p \}$ fixed. It is clear that elements
of $Sh^{(p)}_{(k-r, k)}$ commute with $z_{(p-q+1, p)}$ and
      $(1_{p-q}, -1_q, 1_{k-p})$. Also, $z_{(k-r+1, k)}$ and $(1_{p-q}, -1_q, 1_{k-p})$
commute. Therefore,
\[
\begin{aligned}
    s.s.&( r_{M,G} \circ i_{G,N}(\pi)) \\
     &=
        \sum_{q=0}^{p}      
      \sum_{r=0}^{k-p}
     Sh_{(p-q, p,  k)} Sh^{(p)}_{(k-r, k)}
     z_{(p-q+1, p)}
      z_{(k-r+1, k)}(1_{p-q}, -1_q, 1_{n-p}) (1_{n-r}, -1_r) r_{M,N}(\pi) \\
     &=
        \sum_{q=0}^{p}      
      \sum_{r=0}^{k-p}
      Sh_{(p-q, p, k-r, k)} z_{(p-q+1, p)} z_{(k-r+1, k)}
      (1_{p-q}, -1_q, 1_{k-p-r}, -1_r) 
      r_{M,N}(\pi).
\end{aligned}
\]
\hfill  $\square$

\section{Arthur parameters and Aubert involution: Reducibility at $\frac{1}{2}$}

We study symmetry of the $A$-parameters (\ref{D1})
under the Aubert involution.
Let $\rho$ be an irreducible unitary supercuspidal representation of $GL(n,F)$
and $\sigma$ an irreducible supercuspidal generic representation of $SO(2\ell+1,F)$.
Symmetry of the parameters depends on 
the 
point of reducibility $\alpha$ of the induced representation 
$\nu^\alpha \rho \rtimes \sigma$.
If $\tilde{\rho} \ncong \rho$, then $\nu^\alpha \rho \rtimes \sigma$ is irreducible,
for any $\alpha \in \mathbb{R}$.
Assume $\tilde{\rho} \cong \rho$. 
Then there exists $\alpha \in \{ 0, \frac{1}{2}, 1 \}$ such that 
$\nu^{\pm \alpha} \rho \rtimes \sigma$ is reducible and 
$\nu^\beta \rho \rtimes \sigma$ is irreducible for $|\beta| \ne \alpha$
\cite{Shah:proof}. 
In this section and two consecutive sections, we consider the cases
$\alpha= \frac{1}{2}$, 0 and 1.

\begin{theorem}\label{Th1}
Let $\rho$ be an irreducible unitary supercuspidal representation of $GL(n,F)$
and $\sigma$ an irreducible supercuspidal generic representation of $SO(2\ell+1,F)$.
Suppose  $\nu^\frac{1}{2} \rho \rtimes \sigma$ is reducible.
Let
$
     \bigoplus_{i \in A} \phi_i \otimes S_1
$
be the L-parameter of $\sigma$ and 
$\phi$ be the L-parameter of $\rho$.
Let $\pi$ be the representation with the A-parameter 
\begin{equation}\label{D1}
    {\psi}= \phi \otimes S_{k} \otimes S_{2}
\oplus
     \bigoplus_{i \in A} \phi_i \otimes S_1 \otimes S_1,
\end{equation}
$k \geq 1$. Let $\hat{\psi}$ be the A-parameter of $\hat{\pi}$.

\begin{enumerate}
 
  \item[{(i)}] If $k$ is even, then
\[
  \hat{\psi} =  \phi \otimes S_1 \otimes S_{k+1} \oplus \phi \otimes S_1 \otimes S_{k-1}
\oplus
     \bigoplus_{i \in A} \phi_i \otimes S_1 \otimes S_1.
\]

  \item[{(ii)}] If $k$ is odd, then
\[
   \hat{\psi} = \phi \otimes S_{2} \otimes S_{k} \oplus 
      \bigoplus_{i \in A} \phi_i \otimes S_1 \otimes S_1.
\]

\end{enumerate}

\end{theorem}

{\it Proof.} (i) Let $k=2m$ even.  According to (\ref{B}),
$
    \pi =L_s(\delta [\nu^{-m} \rho ,\nu^{m-1}\rho], \sigma).
$
This means that $\pi$ is the unique irreducible subrepresentation of 
$
   \delta [\nu^{-m} \rho ,\nu^{m-1}\rho]\rtimes \sigma.
$
Theorem 9.1 of \cite{Tad:israel}
tells us 
$
   \delta [\nu^{-m} \rho ,\nu^{m-1}\rho]\rtimes \sigma
$
is irreducible. It follows
\[
   \pi = \delta [\nu^{-m} \rho ,\nu^{m-1}\rho]\rtimes \sigma
\]
and $\pi$ is generic \cite{Muic}. We can apply Lemma 4.2 of 
\cite{BanZh} and Theorem 6.1 of \cite{J-S} to find the $L$-parameter
of $\hat{\pi}$. The $L$-parameter
of $\hat{\pi}$ is equal to
\[
   \phi_{\hat{\psi}} = \bigoplus_{j=-m}^{m-1} (|\cdot|^j \phi \otimes S_1
    \oplus  |\cdot|^{-j}\phi \otimes S_1) \oplus
     \bigoplus_{i \in A} \phi_i \otimes S_1.
\]
Before we apply (\ref{A}), we have to rearrange this expression:
\[
   \phi_{\hat{\psi}} =  \bigoplus_{j=-m}^{m} |\cdot|^j \phi \otimes S_1
    \oplus \bigoplus_{j=-m+1}^{m-1} |\cdot|^j \phi \otimes S_1 \oplus
     \bigoplus_{i \in A} \phi_i \otimes S_1.
\]
Now, (\ref{A}) implies
\[
  \hat{\psi} =  \phi \otimes S_1 \otimes S_{2m+1} \oplus \phi \otimes S_1 \otimes S_{2m-1}
\oplus
     \bigoplus_{i \in A} \phi_i \otimes S_1 \otimes S_1.
\]

(ii) Let $k=2m+1$ odd. The proof is by induction on $m \geq 0$.
For $m=0$, the parameter 
$
\hat{\psi} = \phi \otimes S_{2} \otimes S_{1} \oplus 
      \bigoplus_{i \in A} \phi_i \otimes S_1 \otimes S_1
$
is tempered. It follows from \cite{BanZh} that
$
\widehat{\hat{\psi}} = \phi \otimes S_{1} \otimes S_{2} \oplus 
      \bigoplus_{i \in A} \phi_i \otimes S_1 \otimes S_1.
$

Now, assume the theorem holds for $m-1$ and prove it holds for $m$.
From (\ref{B}),
\[
    \pi =L_s(\delta [\nu^{-m-\frac{1}{2}} \rho ,\nu^{m-\frac{1}{2}}\rho], \sigma).
\]
Let $\tau$ be the representation corresponding to the $A$-parameter
$
        \phi \otimes S_{2} \otimes S_{2m+1} \oplus 
      \bigoplus_{i \in A} \phi_i \otimes S_1 \otimes S_1.
$
Then, by (\ref{A}), the corresponding $L$-parameter is 
\[
 \begin{aligned}
  \bigoplus_{j = -m}^{m} &|\cdot|^j \phi \otimes S_{2} \oplus 
      \bigoplus_{i \in A} \phi_i \otimes S_1  \\
 &= \bigoplus_{j = 1}^{m}( |\cdot|^j \phi \otimes S_{2} \oplus |\cdot|^{-j} \phi \otimes S_{2} )
  \oplus \phi \otimes S_{2} \oplus
      \bigoplus_{i \in A} \phi_i \otimes S_1 
\end{aligned}
\]
and
$
  \tau = L_s( \delta[\nu^{-m-\frac{1}{2}} \rho ,\nu^{-m+\frac{1}{2}}\rho], 
  \delta[\nu^{-m+\frac{1}{2}} \rho ,\nu^{-m+\frac{3}{2}}\rho],
   \cdots,
  \delta[\nu^{-\frac{3}{2}} \rho ,\nu^{-\frac{1}{2}}\rho],
  \delta(\nu^{\frac{1}{2}} \rho; \sigma)).
$
We have to prove $\hat{\pi} = \tau$. Let
\[
  \begin{aligned}
    \pi_1 &=L_s(\delta [\nu^{-m+\frac{1}{2}} \rho ,\nu^{m-\frac{3}{2}}\rho], \sigma), \\
  \tau_1 &= L_s(  
  \delta[\nu^{-m+\frac{1}{2}} \rho ,\nu^{-m+\frac{3}{2}}\rho],
   \cdots,
  \delta[\nu^{-\frac{3}{2}} \rho ,\nu^{-\frac{1}{2}}\rho],
  \delta(\nu^{\frac{1}{2}} \rho; \sigma)).
  \end{aligned}
\]
By the induction assumption, $\hat{\pi}_1 = \tau_1$.
To apply the assumption, we have to prove 
\begin{equation}\label{E}
  \pi \hookrightarrow
   \nu^{m-\frac{1}{2}}\rho \times \nu^{m+\frac{1}{2}}\rho
  \rtimes L_s(\delta [\nu^{-m+\frac{1}{2}} \rho ,
\nu^{m-\frac{3}{2}}\rho], \sigma) =
   \nu^{m-\frac{1}{2}}\rho \times \nu^{m+\frac{1}{2}}\rho \rtimes \pi_1.
\end{equation}
We do it in two steps. First,
define
\[
  \begin{aligned}
   \Pi_1 &= \nu^{m-\frac{1}{2}}\rho \rtimes L_s(\delta [\nu^{-m-\frac{1}{2}} \rho ,
\nu^{m-\frac{3}{2}}\rho], \sigma), \quad
  \Pi = \nu^{m-\frac{1}{2}}\rho \times 
   \delta [\nu^{-m-\frac{1}{2}} \rho ,\nu^{m-\frac{3}{2}}\rho] \rtimes \sigma \quad 
  \end{aligned}
\]
and
\[
   \Pi_0 = \nu^{m-\frac{1}{2}}\rho \otimes \nu^{m-\frac{3}{2}}\rho \otimes \cdots
   \otimes \nu^{-m+\frac{1}{2}} \rho \otimes \nu^{-m-\frac{1}{2}} \rho \otimes \sigma.
\]
Denote by $M$ the standard Levi subgroup of $G= SO(2(n(2m+1)+\ell)+1, F)$ isomorphic to 
$
     GL(n,F)^{2m+1} \times SO(2\ell+1, F).
$
Let
$
    N \cong GL(n,F) \times GL(2mn,F) \times SO(2\ell+1, F)
$
be an intermediate Levi subgroup, $M<N<G$. 
Then 
\[
    r_{M,N}(\nu^{m-\frac{1}{2}}\rho \otimes 
   \delta [\nu^{-m-\frac{1}{2}} \rho ,\nu^{m-\frac{3}{2}}\rho] \otimes \sigma)
   = \Pi_0
\]
and,  according to Lemma \ref{Can2},
\begin{equation}\label{C}
\begin{aligned}
    s.s.&( r_{M,G} (\Pi)) \\
     &=
        \sum_{q=0}^{1}      
      \sum_{r=0}^{2m}
      Sh_{(1-q, 1, 2m+1-r, 2m+1)}  z_{(2m-r+2, 2m+1)}
      (1_{1-q}, -1_q, 1_{2m-r}, -1_r) 
      \Pi_0.
\end{aligned}
\end{equation}
Here, we use the fact that $z_{(2-q, 1)} =1$ for $q=0$ or $q=1$.
We claim the multiplicity of $\Pi_0$ in $r_{M,G} (\Pi)$ is one.
Indeed, $\Pi_0$ appears in (\ref{C}) only for $r=0$.
This comes from the fact that $\nu^{m+\frac{1}{2}}\rho$ does not appear in $\Pi_0$.
Therefore, we have to consider
$ 
       \sum_{q=0}^{1} 
      Sh_{(1, 2m+1)}  
      (1_{1-q}, -1_q, 1_{2m-r}) 
      \Pi_0.
$
If $q=1$,  then the first factor $\nu^{m-\frac{1}{2}}\rho$ changes into
$\nu^{-m+\frac{1}{2}}\rho$
and
$
      Sh_{(1, 2m+1)}  
      \nu^{-m+\frac{1}{2}}\rho \otimes \nu^{m-\frac{3}{2}}\rho \otimes \cdots
   \otimes \nu^{-m+\frac{1}{2}} \rho \otimes \nu^{-m-\frac{1}{2}} \rho \otimes \sigma
$
cannot produce $\Pi_0$.
If $q=0$,   then
$      Sh_{(1, 2m+1)}  
       \Pi_0
$
produces $\Pi_0$ only for $1 \in Sh_{(1, 2m+1)}$.

Note that $\pi$ and $\Pi_1$ are subrepresentations of $\Pi$.
We see easily that $\Pi_0 \leq r_{M,G}(\pi)$ and $\Pi_0 \leq r_{M,G}(\Pi_1)$.
This implies $\pi$ and $\Pi_1$ have a subquotient in common.
In other words (since $\pi$ is irreducible), $\pi$ is a subquotient of $\Pi_1$.
Lemma \ref{Can3} implies $\pi$ is a subrepresentation of $\Pi_1$.

Next, we consider $\pi' = L_s(\delta [\nu^{-m+\frac{1}{2}} \rho ,
\nu^{m-\frac{3}{2}}\rho], \sigma)$ and define
\[
  \begin{aligned}
   \Pi_1' &= \nu^{m+\frac{1}{2}}\rho \rtimes L_s(\delta [\nu^{-m+\frac{1}{2}} \rho ,
\nu^{m-\frac{3}{2}}\rho], \sigma), \quad
  \Pi' = \nu^{m+\frac{1}{2}}\rho \times 
   \delta [\nu^{-m+\frac{1}{2}} \rho ,\nu^{m-\frac{3}{2}}\rho] \rtimes \sigma,
  \end{aligned}
\]
and
\[
   \Pi_0' = \nu^{m+\frac{1}{2}}\rho \otimes (\nu^{m-\frac{3}{2}}\rho \otimes \cdots
   \otimes \nu^{-m+\frac{1}{2}} \rho ) \otimes \sigma.
\]
Denote by $M'$ the standard Levi subgroup of $G'= SO(2(2mn+\ell)+1, F)$ isomorphic to 
$
     GL(n,F)^{2m} \times SO(2\ell+1, F).
$
Let
$
    N' \cong GL(n,F) \times GL((2m-1)n,F) \times SO(2\ell+1, F)
$
be an intermediate Levi subgroup, $M'<N'<G'$. 
Then 
\[
    r_{M',N'}(\nu^{m+\frac{1}{2}}\rho \otimes 
   \delta [\nu^{-m+\frac{1}{2}} \rho ,\nu^{m-\frac{3}{2}}\rho] \otimes \sigma)
   = \Pi_0'
\]
and,  according to Lemma \ref{Can2},
\begin{equation}\label{D}
\begin{aligned}
    s.s.&( r_{M',G'} (\Pi')) \\
     &=
        \sum_{q=0}^{1}      
      \sum_{r=0}^{2m-1}
      Sh_{(1-q, 1, 2m-r, 2m)}  z_{(2m-r+1, 2m)}
      (1_{1-q}, -1_q, 1_{2m-r-1}, -1_r) 
      \Pi_0'.
\end{aligned}
\end{equation}
Since $\nu^{-m-\frac{1}{2}}\rho$ is not present in $\Pi_0'$,
it is obvious that $\Pi_0'$ appears in (\ref{D}) only for $q=0$.
Further,  $\Pi_0'$ appears in
$     \sum_{r=0}^{2m-1}
      Sh_{(1, 2m-r, 2m)}  z_{(2m-r+1, 2m)}
      (1_{2m-r}, -1_r) 
      \Pi_0'
$
only for $r=0$. It follows the multiplicity of $\Pi_0'$ in $r_{M',G'} (\Pi')$ is one.
Clearly, 
$
   r_{M',G'} (\pi') \geq 
 (\nu^{m-\frac{3}{2}}\rho \otimes \cdots
   \otimes \nu^{-m+\frac{1}{2}} \rho ) \otimes \nu^{-m-\frac{1}{2}}\rho \otimes \sigma.
$
We apply Lemma \ref{lemma} on $\nu^{-m-\frac{1}{2}}\rho$ to show that
$
   r_{M',G'} (\pi') \geq 
 (\nu^{m-\frac{3}{2}}\rho \otimes \cdots
   \otimes \nu^{-m+\frac{1}{2}} \rho ) \otimes \nu^{m+\frac{1}{2}}\rho \otimes \sigma
$
and 
$
   r_{M',G'} (\pi') \geq \Pi_0'.
$
In a similar way as earlier, we show  $\pi'$ is a subrepresentation of
$\Pi_1'$. 
We have proved 
$
   \pi \hookrightarrow \nu^{m-\frac{1}{2}}\rho \rtimes \pi'
    \hookrightarrow \nu^{m-\frac{1}{2}}\rho \rtimes \nu^{m+\frac{1}{2}}\rho \rtimes \pi_1.
$
This implies (\ref{E}).
Let us mention that the arguments presented here do not work if we try to put
two steps of the proof into one single step, because the multiplicity of
$
  \nu^{m-\frac{1}{2}}\rho \otimes
\nu^{m+\frac{1}{2}}\rho \otimes (\nu^{m-\frac{3}{2}}\rho \otimes \cdots
   \otimes \nu^{-m+\frac{1}{2}} \rho ) \otimes \sigma
$
in the Jacquet module of
$
   \nu^{m-\frac{1}{2}}\rho \times
  \nu^{m+\frac{1}{2}}\rho \times 
   \delta [\nu^{-m+\frac{1}{2}} \rho ,\nu^{m-\frac{3}{2}}\rho] \rtimes \sigma
$
is greater than 2.

The Aubert involution is defined on the Grothendieck group.
It commutes with parabolic induction.
If we apply the 
Aubert involution on (\ref{E}), we see that $\hat{\pi}$ is a component 
of the representation
\begin{equation}\label{F}
   \nu^{m-\frac{1}{2}}\rho \times \nu^{m+\frac{1}{2}}\rho \rtimes \hat{\pi}_1 =
   \nu^{m-\frac{1}{2}}\rho \times \nu^{m+\frac{1}{2}}\rho \rtimes \tau_1.
\end{equation}
We show that  $\tilde{\pi} \cong \pi$ and $\tilde{\tau}_1 \cong \tau_1$.
The representation $\sigma$ is generic and supercuspidal.
Then $\tilde{\sigma}$ is also  generic and supercuspidal.
The representations $\sigma$ and $\tilde{\sigma}$ belong to the same 
$L$-packet. According to Theorem 1.1 of \cite{J-S}, there is a bijection
between the set of equivalence classes of irreducible supercuspidal generic
representations of $SO(2n+1,F)$ and the set of $L$-parameters described in 
Theorem 1.1 of \cite{J-S}. This implies $\tilde{\sigma} \cong \sigma$.
As explained in Section 6 of \cite{Tad:compositio},
$L_s(\delta_1, \dots, \delta_n, \sigma)\tilde{} = 
L_s(\delta_1, \dots, \delta_n, \tilde{\sigma})$.
It follows
$
    \tilde{\pi} =
  L_s(\delta [\nu^{-m-\frac{1}{2}} \rho ,\nu^{m-\frac{1}{2}}\rho], \tilde{\sigma})
\cong 
 L_s(\delta [\nu^{-m-\frac{1}{2}} \rho ,\nu^{m-\frac{1}{2}}\rho], {\sigma})
 = {\pi}.
$
Similarly, $\tilde{\tau_1} \cong \tau_1$.

Therefore, $(\hat{\pi})\tilde{} \cong (\tilde{\pi}) \hat{} \cong \hat{\pi}$.
It follows that $\hat{\pi}$ is a component of the contragradient 
of (\ref{F}), that is, $\hat{\pi}$ is a component of
$
     \nu^{-m+\frac{1}{2}}\rho \times \nu^{-m-\frac{1}{2}}\rho \rtimes \tau_1.
$
On the other hand, Frobenius reciprocity and (\ref{E}) imply
$
    r_{N,G}(\pi) \geq \nu^{m-\frac{1}{2}}\rho \otimes \nu^{m+\frac{1}{2}}\rho
  \otimes \pi_1,
$
where
$
    N \cong GL(n,F) \times GL(n,F) \times SO(2((2m-1)n+\ell)+1, F).
$
From the exactness of the Jacquet functor, we have
\begin{equation}\label{O}
    r_{M,G}(\pi) \geq \nu^{m-\frac{1}{2}}\rho \otimes \nu^{m+\frac{1}{2}}\rho
  \otimes r_{M_2,G_2} (\pi_1),
\end{equation}
where $G_2 = SO(2((2m-1)n+\ell)+1, F)$ and 
$M_2 \cong GL(n,F)^{2m-1} \times SO(2(\ell)+1, F)$.
We apply the Aubert involution on (\ref{O}).
According to \cite{Aub}, Th\'eorem\`e 1.7, 
\[
    r_{M,G}(\hat{\pi}) \geq \nu^{-m+\frac{1}{2}}\rho \otimes \nu^{-m-\frac{1}{2}}\rho
  \otimes r_{M_2,G_2} (\hat{\pi}_1) =
  \nu^{-m+\frac{1}{2}}\rho \otimes \nu^{-m-\frac{1}{2}}\rho
  \otimes r_{M_2,G_2} (\tau_1).
\]
In particular, $r_{M,G}(\hat{\pi}) \geq \Pi_0''$, where
\[
 \Pi_0'' = \nu^{-m+\frac{1}{2}}\rho \otimes \nu^{-m-\frac{1}{2}}\rho
  \otimes 
   (\nu^{-m+\frac{3}{2}} \rho \otimes \nu^{-m+\frac{1}{2}}\rho) \otimes
   \cdots \otimes 
  (\nu^{-\frac{1}{2}} \rho \otimes \nu^{-\frac{3}{2}}\rho)\otimes
  \nu^{\frac{1}{2}} \rho \otimes \sigma.
\]
To finish the proof, we need the following two lemmas.

\begin{lemma}\label{2}
 The representation $\Pi_0''$ appears with multiplicity one in the Jacquet module of
\[
 \Pi'' =
\delta[\nu^{-m-\frac{1}{2}} \rho ,\nu^{-m+\frac{1}{2}}\rho] \times 
  \delta[\nu^{-m+\frac{1}{2}} \rho ,\nu^{-m+\frac{3}{2}}\rho] \times
   \cdots \times
  \delta[\nu^{-\frac{3}{2}} \rho ,\nu^{-\frac{1}{2}}\rho] \times
  \delta(\nu^{\frac{1}{2}} \rho; \sigma).
\]

\end{lemma}

{\it Proof.} By induction on $m \geq 1$. This proof is independent of the proof of  
Theorem \ref{Th1}, and two inductions do not interfere. 

Let 
$
    N \cong GL(2n,F) \times SO(2((2m-1)n+\ell)+1, F)
$
be an intermediate Levi subgroup, $M < N < G$.
Assume $m=1$. Then, by Lemma \ref{Can},
\begin{equation}\label{H}
\begin{aligned}
    s.s.&( r_{M,G} (\Pi'')) \\
      &=
        \sum_{q=0}^{2}  
      Sh_{(2-q, 2, 3)} z_{(3-q, 2)} (1_{2-q}, -1_q, 1)\,
      \nu^{-\frac{1}{2}} \rho \otimes \nu^{-\frac{3}{2}}\rho \otimes
  \nu^{\frac{1}{2}} \rho \otimes \sigma.
\end{aligned}
\end{equation}
The representation 
$
\nu^{-\frac{1}{2}} \rho \otimes \nu^{-\frac{3}{2}}\rho \otimes
  \nu^{\frac{1}{2}} \rho \otimes \sigma
$
appears in (\ref{H}) only for $q=0$, for the permutation $1 \in Sh_{(2, 2m+1)}$.
Therefore, the multiplicity of 
$
\nu^{-\frac{1}{2}} \rho \otimes \nu^{-\frac{3}{2}}\rho \otimes
  \nu^{\frac{1}{2}} \rho \otimes \sigma
$
in $r_{M,G} (\Pi'')$ is one.

Now, assume the lemma holds for $m-1$. 
From Lemma \ref{Can},
\begin{equation}\label{G}
  \begin{aligned}
    s.s.&( r_{M,G} (\Pi'')) =
        \sum_{q=0}^{2}  
      Sh_{(2-q, 2, 2m+1)} z_{(3-q, 2)} (1_{2-q}, -1_q, 1_{2m-1})\,\\
      &\nu^{-m+\frac{1}{2}} \rho \otimes \nu^{-m-\frac{1}{2}}\rho 
  \otimes 
 r_{M',G'}( \delta[\nu^{-m+\frac{1}{2}} \rho ,\nu^{-m+\frac{3}{2}}\rho] \times
   \cdots \times
  \delta[\nu^{-\frac{3}{2}} \rho ,\nu^{-\frac{1}{2}}\rho] \times
  \delta(\nu^{\frac{1}{2}} \rho; \sigma)).
\end{aligned}
\end{equation}
Observe that neither $\nu^{-m-\frac{1}{2}}\rho$ nor $\nu^{m+\frac{1}{2}}\rho$ appear in 
\begin{equation}\label{K}
 r_{M',G'}( \delta[\nu^{-m+\frac{1}{2}} \rho ,\nu^{-m+\frac{3}{2}}\rho] \times
   \cdots \times
  \delta[\nu^{-\frac{3}{2}} \rho ,\nu^{-\frac{1}{2}}\rho] \times
  \delta(\nu^{\frac{1}{2}} \rho; \sigma)).
\end{equation}
Therefore, $\Pi''_0$ appears in (\ref{G}) only for $q=0$ and for all 
permutations $s \in Sh_{(2, 2m+1)}$ which keep $\nu^{-m-\frac{1}{2}}\rho$
fixed. There is only one such permutation, namely $s =1$.
It follows  the multiplicity of $\Pi''_0$ in $r_{M,G} (\Pi'')$
is equal to the multiplicity of 
\[ 
   (\nu^{-m+\frac{3}{2}} \rho \otimes \nu^{-m+\frac{1}{2}}\rho) \otimes
   \cdots \otimes 
  (\nu^{-\frac{1}{2}} \rho \otimes \nu^{-\frac{3}{2}}\rho)\otimes
  \nu^{\frac{1}{2}} \rho \otimes \sigma
\]
in (\ref{K}) which is, by the induction assumption, equal to one.
\hfill $\square$

\begin{lemma}\label{3}
  The representation $\Pi_0''$ appears with multiplicity one in the Jacquet module of
\[
 \Pi''' =
\nu^{-m+\frac{1}{2}} \rho \times \nu^{-m-\frac{1}{2}}\rho \times 
  \delta[\nu^{-m+\frac{1}{2}} \rho ,\nu^{-m+\frac{3}{2}}\rho] \times
   \cdots \times
  \delta[\nu^{-\frac{3}{2}} \rho ,\nu^{-\frac{1}{2}}\rho] \times
  \delta(\nu^{\frac{1}{2}} \rho; \sigma).
\]

\end{lemma}

{\it Proof.}  Let 
$
    N \cong GL(n,F) \times GL(n,F) \times SO(2((2m-1)n+\ell)+1, F)
$
be an intermediate Levi subgroup, $M < N < G$.
Straightforward computation shows the lemma holds for $m=1$.
Now, assume $m >1$.
 Then,
\begin{equation}\label{L}
\begin{aligned}
    s.s.&( r_{M,G} (\Pi''')) 
     =
        \sum_{q=0}^{1}      
      \sum_{r=0}^{1}
      Sh_{(1,2, 2m+1)} 
      (1_{1-q}, -1_q, 1_{1-r}, -1_r, 1_{2m-1}) \\
      &\nu^{-m+\frac{1}{2}} \rho \otimes \nu^{-m-\frac{1}{2}}\rho \otimes 
 r_{M',G'}( \delta[\nu^{-m+\frac{1}{2}} \rho ,\nu^{-m+\frac{3}{2}}\rho] \times
   \cdots \times
  \delta[\nu^{-\frac{3}{2}} \rho ,\nu^{-\frac{1}{2}}\rho] \times
  \delta(\nu^{\frac{1}{2}} \rho; \sigma)).
\end{aligned}
\end{equation}
To obtain $\Pi_0''$ in (\ref{L}), we obviously need $q=0$ and $r=0$.
Suppose 
\begin{equation}\label{M}
\Pi_0'' = s(\nu^{-m+\frac{1}{2}} \rho \otimes \nu^{-m-\frac{1}{2}}\rho \otimes 
    \chi),
\end{equation}
where $ s \in Sh_{(1,2, 2m+1)}$ and  $\chi \leq 
r_{M',G'}( \delta[\nu^{-m+\frac{1}{2}} \rho ,\nu^{-m+\frac{3}{2}}\rho] \times
   \cdots \times
  \delta[\nu^{-\frac{3}{2}} \rho ,\nu^{-\frac{1}{2}}\rho] \times
  \delta(\nu^{\frac{1}{2}} \rho; \sigma))$.
If $s \ne 1$, then (\ref{M}) is possible only for
\[
   \chi =   (\nu^{-m+\frac{1}{2}} \rho \otimes \nu^{-m+\frac{3}{2}}\rho) \otimes
   \cdots \otimes 
  (\nu^{-\frac{1}{2}} \rho \otimes \nu^{-\frac{3}{2}}\rho)\otimes
  \nu^{\frac{1}{2}} \rho \otimes \sigma.
\]
This representation, however, does not appear in 
\begin{equation}\label{N}
r_{M',G'}( \delta[\nu^{-m+\frac{1}{2}} \rho ,\nu^{-m+\frac{3}{2}}\rho] \times
   \cdots \times
  \delta[\nu^{-\frac{3}{2}} \rho ,\nu^{-\frac{1}{2}}\rho] \times
  \delta(\nu^{\frac{1}{2}} \rho; \sigma)).
\end{equation}
It follows the multiplicity of $\Pi_0''$ in $r_{M,G} (\Pi''')$ is equal to the
multiplicity of 
$
(\nu^{-m+\frac{3}{2}} \rho \otimes \nu^{-m+\frac{1}{2}}\rho) \otimes
   \cdots \otimes 
  (\nu^{-\frac{1}{2}} \rho \otimes \nu^{-\frac{3}{2}}\rho)\otimes
  \nu^{\frac{1}{2}} \rho \otimes \sigma
$
in (\ref{N}) which is, by Lemma \ref{2}, equal to one.   \hfill $\square$

To complete the proof of Theorem \ref{Th1}, observe first that
Lemmas \ref{2} and \ref{3} imply $\hat{\pi}$ is a subrepresentation of 
$\Pi''$. Since $\Pi''$ is the representation induced from Langlands data,
it has a unique subrepresentation. Therefore,
\[
   \hat{\pi} =
 L_s( \delta[\nu^{-m-\frac{1}{2}} \rho ,\nu^{-m+\frac{1}{2}}\rho], 
  \delta[\nu^{-m+\frac{1}{2}} \rho ,\nu^{-m+\frac{3}{2}}\rho],
   \cdots,
  \delta[\nu^{-\frac{3}{2}} \rho ,\nu^{-\frac{1}{2}}\rho],
  \delta(\nu^{\frac{1}{2}} \rho; \sigma)) = \tau,
\]
finishing the proof.
\hfill  $\square$

Directly from the proof of the theorem, we have the following:

\begin{corollary}
Let $\rho$ be an irreducible unitary supercuspidal representation of $GL(n,F)$
and $\sigma$ an irreducible supercuspidal generic representation of $SO(2\ell+1,F)$.
Suppose  $\nu^\frac{1}{2} \rho \rtimes \sigma$ is reducible.
 \begin{enumerate}
     \item[{(i)}] If 
$
    \pi =L_s(\delta [\nu^{-m} \rho ,\nu^{m-1}\rho], \sigma),
$
then the Aubert involution of $\pi$ is equal to
\[
   \hat{\pi} 
= L_s(\nu^{-m} \rho, \nu^{-m+1} \rho,\nu^{-m+1} \rho, \dots, \nu^{-1} \rho,\nu^{-1} \rho, 
    \rho \rtimes \sigma). 
\]
\item[{(ii)}] If
$
    \pi =L_s(\delta [\nu^{-m-\frac{1}{2}} \rho ,\nu^{m-\frac{1}{2}}\rho], \sigma),
$
then
\[
  \hat{\pi} = L_s( \delta[\nu^{-m-\frac{1}{2}} \rho ,\nu^{-m+\frac{1}{2}}\rho], 
  \delta[\nu^{-m+\frac{1}{2}} \rho ,\nu^{-m+\frac{3}{2}}\rho],
   \cdots,
  \delta[\nu^{-\frac{3}{2}} \rho ,\nu^{-\frac{1}{2}}\rho],
  \delta(\nu^{\frac{1}{2}} \rho; \sigma)).
\]

\end{enumerate}

\end{corollary}

\section{Arthur parameters and Aubert involution: Reducibility at 0}

We continue to study symmetry of the $A$-parameters 
under the Aubert involution.

\begin{theorem}\label{Th2}
Let $\rho$ be an irreducible unitary supercuspidal representation of $GL(n,F)$
and $\sigma$ an irreducible supercuspidal generic representation of $SO(2\ell+1,F)$.
Suppose  $\rho \rtimes \sigma$ is reducible.
Let
$
     \bigoplus_{i \in A} \phi_i \otimes S_1
$
be the L-parameter of $\sigma$ and 
$\phi$ be the L-parameter of $\rho$.
Let $\pi$ be the representation with the A-parameter 
\[
    {\psi}= \phi \otimes S_{k} \otimes S_{2}
\oplus
     \bigoplus_{i \in A} \phi_i \otimes S_1 \otimes S_1,
\]
$k \geq 1$. Let $\hat{\psi}$ be the A-parameter of $\hat{\pi}$.

\begin{enumerate}
 
\item[{(i)}] Assume $k$ is odd. If $k=1$, then $\hat{\pi} = \pi$.
If $k \geq 3$, then
\[
  \hat{\psi} =  \phi \otimes S_1 \otimes S_{k+1} \oplus \phi \otimes S_1 \otimes S_{k-1}
\oplus
     \bigoplus_{i \in A} \phi_i \otimes S_1 \otimes S_1.
\]

\item[{(ii)}] If $k$ is even, then
\[
   \hat{\psi} = \phi \otimes S_{2} \otimes S_{k} \oplus 
      \bigoplus_{i \in A} \phi_i \otimes S_1 \otimes S_1.
\]

\end{enumerate}

\end{theorem}

{\it Proof.} (i) Let $k=2m+1$ odd. The proof is similar to the proof of Theorem \ref{Th1}, (i).
According to (\ref{B}), $\pi$ is equal to
$
    L_s(\delta [\nu^{-m-\frac{1}{2}} \rho ,\nu^{m-\frac{1}{2}}\rho], \sigma).
$
Theorem 9.1 of \cite{Tad:israel}
tells us 
$
   \delta [\nu^{-m-\frac{1}{2}} \rho ,\nu^{m-\frac{1}{2}}\rho]\rtimes \sigma
$
is irreducible. It follows
$
   \pi = \delta [\nu^{-m-\frac{1}{2}} \rho ,\nu^{m-\frac{1}{2}}\rho] \rtimes \sigma
$
and $\pi$ is generic \cite{Muic}.  
If $m=0$, then $\pi = \nu^{\frac{1}{2}} \rtimes \sigma$
and $\hat{\pi} = \pi$.
Assume $m \geq 1$. Then the $L$-parameter
of $\hat{\pi}$ is equal to
\[
 \begin{aligned}
   \phi_{\hat{\psi}} &= \bigoplus_{j=-m-\frac{1}{2}}^{m-\frac{1}{2}} 
(|\cdot|^j \phi \otimes S_1
    \oplus  |\cdot|^{-j} \phi \otimes S_1) \oplus
     \bigoplus_{i \in A} \phi_i \otimes S_1 \\
   &=  \bigoplus_{j=-m-\frac{1}{2}}^{m+\frac{1}{2}} |\cdot|^j \phi \otimes S_1
    \oplus \bigoplus_{j=-m+\frac{1}{2}}^{m-\frac{1}{2}} |\cdot|^j\phi \otimes S_1 \oplus
     \bigoplus_{i \in A} \phi_i \otimes S_1.
   \end{aligned}
\]
Now, (\ref{A}) implies
$
  \hat{\psi} =  \phi \otimes S_1 \otimes S_{2m+2} \oplus \phi \otimes S_1 \otimes S_{2m}
\oplus
     \bigoplus_{i \in A} \phi_i \otimes S_1 \otimes S_1.
$

(ii) Let $k=2m$ even. The proof is by induction on $ m\geq 1$. Assume $m=1$. Then
$
   \psi = \phi \otimes S_{2} \otimes S_2 \oplus 
      \bigoplus_{i \in A} \phi_i \otimes S_1 \otimes S_1
$
and we have to prove $\pi = \hat{\pi}$.
According to (\ref{B}), 
$\pi=L_s(\delta[\nu^{-1}\rho, \rho], \sigma)$. Now,
\[
    r_{M,G}(\delta[\nu^{-1}\rho, \rho] \rtimes \sigma) = 
\rho \otimes \nu^{-1}\rho \otimes \sigma + \rho \otimes \nu\rho \otimes \sigma
     + 2 \nu \rho \otimes \rho \otimes \sigma,
\]
from Lemma \ref{Can}.
Proposition 4.2 of \cite{Tad:compositio} tells us 
$\delta[\nu^{-1}\rho, \rho]\rtimes \sigma$
and 
$
\delta[\nu^{-1}\rho, \rho]\,\tilde{} \rtimes \sigma = 
\delta[\rho, \nu \rho] \rtimes \sigma 
$ 
have the same irreducible
components.
The representation $\delta[\rho, \nu \rho] \rtimes \sigma$ is reducible.
It has two discrete series subrepresentations $\tau_1$, $\tau_2$
and the unique Langlands quotient $\pi = L_q(\delta[\rho, \nu \rho] , \sigma)$.
Since 
$
 \pi \hookrightarrow \delta[\nu^{-1}\rho, \rho] \rtimes \sigma,
$
it follows $r_{M,G}(\pi) \geq \rho \otimes \nu^{-1}\rho \otimes \sigma $.
Lemma \ref{lemma} implies 
$r_{M,G}(\pi) \geq \rho \otimes \nu \rho \otimes \sigma $.
We conclude
$
   \delta[\rho, \nu \rho] \rtimes \sigma =  \tau_1 + \tau_2 + \pi
$
and $r_{M,G}(\tau_1) = r_{M,G}(\tau_2) = \nu \rho \otimes \rho \otimes \sigma$,
\begin{equation}\label{P}
   r_{M,G}(\pi) = \rho \otimes \nu^{-1}\rho \otimes \sigma +
    \rho \otimes \nu \rho \otimes \sigma.
\end{equation}
We apply the Aubert involution on (\ref{P}).
According to \cite{Aub}, Th\'eorem\`e 1.7, 
\begin{equation}\label{Q}
   r_{M,G}(\hat{\pi}) = \rho \otimes \nu \rho \otimes \sigma +
      \rho \otimes \nu^{-1}\rho \otimes \sigma.
\end{equation}
We search in (\ref{Q}) for representations coming from Langlands data
in subrepresentation setting. The representation
$
\rho \otimes \nu \rho \otimes \sigma
$
does not come from Langlands data.
The representation 
$
\rho \otimes \nu^{-1}\rho \otimes \sigma
$
comes from the Langlands data $\delta[\nu^{-1}\rho, \rho]\otimes \sigma$,
which is precisely the Langlands data for $\pi$. 
It follows $\pi= \hat{\pi}$.

Now, assume the theorem holds for $m$ and prove it holds for $m+1$.
This case is similar to Theorem \ref{Th1}, (ii). For that reason, 
we skip the detail and give an outline of the proof.
From (\ref{B}),
$
    \pi =L_s(\delta [\nu^{-m-1} \rho ,\nu^{m}\rho], \sigma).
$
Let $\tau$ be the representation corresponding to the $A$-parameter
$
        \phi \otimes S_{2} \otimes S_{2m+2} \oplus 
      \bigoplus_{i \in A} \phi_i \otimes S_1 \otimes S_1.
$
Then, by (\ref{A}), the corresponding $L$-parameter is 
\[
 \begin{aligned}
  \bigoplus_{j = -m-\frac{1}{2}}^{m+\frac{1}{2}} &|\cdot|^j \phi \otimes S_{2} \oplus 
      \bigoplus_{i \in A} \phi_i \otimes S_1  
 = \bigoplus_{j = \frac{1}{2}}^{m+\frac{1}{2}}
  ( |\cdot|^j \phi \otimes S_{2} \oplus |\cdot|^{-j} \phi \otimes S_{2} )
  \oplus
      \bigoplus_{i \in A} \phi_i \otimes S_1 
\end{aligned}
\]
and
$
  \tau = L_s( \delta[\nu^{-m-1} \rho ,\nu^{-m}\rho], 
  \delta[\nu^{-m} \rho ,\nu^{-m+1}\rho],
   \dots,
  \delta[\nu^{-1} \rho ,\rho],
  \sigma).
$
We have to prove $\hat{\pi} = \tau$. 
It can be shown that
$
   \Pi_0 = \nu^{m}\rho \otimes \nu^{m-1}\rho \otimes \cdots
   \otimes \nu^{-m} \rho \otimes \nu^{-m-1} \rho \otimes \sigma
$
appears with multiplicity one in the Jacquet modules of each of the following representations:
\[
   \pi, \quad
   \Pi_1 = \nu^{m}\rho \rtimes L_s(\delta [\nu^{-m-1} \rho ,
\nu^{m-1}\rho], \sigma), \quad 
  \Pi = \nu^{m}\rho \times 
   \delta [\nu^{-m-1} \rho ,\nu^{m-1}\rho] \rtimes \sigma. 
\]
This implies 
$
   \pi \hookrightarrow \nu^{m}\rho \rtimes L_s(\delta [\nu^{-m-1} \rho ,
\nu^{m-1}\rho], \sigma).
$

Next, we consider $\pi' = L_s(\delta [\nu^{-m-1} \rho ,
\nu^{m-1} \rho], \sigma)$ and define
\[
   \Pi_1' = \nu^{m+1}\rho \rtimes L_s(\delta [\nu^{-m} \rho ,
\nu^{m-1}], \sigma), \quad
  \Pi' = \nu^{m+1}\rho \times 
   \delta [\nu^{-m} \rho ,\nu^{m-1}\rho] \rtimes \sigma.
\]
It can be shown that
$
   \Pi_0' = \nu^{m+1}\rho \otimes (\nu^{m-1}\rho \otimes \cdots
   \otimes \nu^{-m} \rho ) \otimes \sigma
$
appears with multiplicity one in the Jacquet modules of each of the  representations
$\pi'$, $\Pi_1'$ and $\Pi'$. 
This implies 
$
   \pi' \hookrightarrow \Pi_1'
$
and
\begin{equation}\label{R}
   \pi  \hookrightarrow \nu^{m}\rho \times \nu^{m+1}\rho
\rtimes L_s(\delta [\nu^{-m} \rho ,
\nu^{m-1}\rho], \sigma). 
\end{equation}
If we apply the 
Aubert involution and contragredient
on (\ref{R}), we conclude that $\hat{\pi}$ is a component 
of the representation
$
     \nu^{-m}\rho \times \nu^{-m-1}\rho \rtimes 
   L_s( 
  \delta[\nu^{-m} \rho ,\nu^{-m+1}\rho],
   \dots,
  \delta[\nu^{-1} \rho ,\rho],
  \sigma).
$
On the other hand,  (\ref{R}) implies
$
    r_{N,G}(\pi) \geq \nu^{m}\rho \otimes \nu^{m+1}\rho
  \otimes L_s(\delta [\nu^{-m} \rho ,\nu^{m-1}\rho], \sigma)$
and
\[
    r_{M,G}(\hat{\pi}) \geq \nu^{-m}\rho \otimes \nu^{-m-1}\rho
  \otimes 
   (\nu^{-m+1} \rho \otimes \nu^{-m}\rho) \otimes
   \cdots \otimes 
  ( \rho \otimes \nu^{-1}\rho)\otimes
   \sigma.
\]
Based on considerations of Jacquet modules, we can show
$\hat{\pi} \hookrightarrow \delta[\nu^{-m-1} \rho ,\nu^{-m}\rho] \times 
  \delta[\nu^{-m} \rho ,\nu^{-m+1}\rho] \times 
   \cdots \times
  \delta[\nu^{-1} \rho ,\rho] \rtimes
  \sigma,
$
which proves $\hat{\pi} = \tau$.   \hfill  $\square$

\section{Arthur parameters and Aubert involution: Reducibility at 1}

In this section, we study symmetry of the $A$-parameters 
for the remaining case, for the point of reducibility $\alpha = 1$.
At the end of the section, we give an example where the Aubert involution
sends a base point to a representation which is not a base point (Remark 6.1).

\begin{theorem}\label{Th3}
Let $\rho$ be an irreducible unitary supercuspidal representation of $GL(n,F)$
and $\sigma$ an irreducible supercuspidal generic representation of $SO(2\ell+1,F)$.
Suppose  $\nu \rho \rtimes \sigma$ is reducible.
Let
$
     \bigoplus_{i \in A} \phi_i \otimes S_1
$
be the L-parameter of $\sigma$ and 
$\phi$ be the L-parameter of $\rho$.
Let $\pi$ be the representation with the A-parameter 
\[
    {\psi}= \phi \otimes S_{k} \otimes S_{2}
\oplus
     \bigoplus_{i \in A} \phi_i \otimes S_1 \otimes S_1, \quad k \geq 1.
\]

\begin{enumerate}
 
\item[{(i)}] If $k=1$, then $\hat{\pi}=\pi$.

\item[{(ii)}] Assume $k\geq 3$ is odd.
Let $\hat{\psi}$ be the A-parameter of $\hat{\pi}$.
Then
\[
  \hat{\psi} =  \phi \otimes S_1 \otimes S_{k+1} \oplus \phi \otimes S_1 \otimes S_{k-1}
\oplus
     \bigoplus_{i \in A} \phi_i \otimes S_1 \otimes S_1.
\]

\item[{(iii)}] Assume $k=2m$ is even.
Let $\hat{\varphi}$ be the $L$-parameter of $\hat{\pi}$.
Then
\[
   \hat{\varphi} = (\bigoplus_{j = \frac{3}{2}}^{m-\frac{1}{2}} 
       |\cdot|^j  \phi \otimes S_2 \oplus |\cdot|^{-j}  \phi \otimes S_2)
     \oplus  \phi \otimes S_3 \oplus  \phi \otimes S_1 \oplus
      \bigoplus_{i \in A} \phi_i \otimes S_1.
\]
In particular, if $m \geq 2$, $\hat{\varphi}$ is not the image of an $A$-parameter, i.e.,
$\hat{\varphi}$ is not of the form $\phi_{\hat{\psi}}$, for an $A$-parameter
$\hat{\psi}$.

\end{enumerate}

\end{theorem}

{\it Proof.} (i), (ii) Identical to the proof of Theorem \ref{Th2}, (i).

(iii)
Assume $m=1$. Then 
$\pi=L_s(\delta[\nu^{-1}\rho, \rho], \sigma)$. We have
\begin{equation}\label{S}
    r_{M,G}(\delta[\nu^{-1}\rho, \rho] \rtimes \sigma) =
\rho \otimes \nu^{-1}\rho \otimes \sigma + \rho \otimes \nu\rho \otimes \sigma
     + 2 \, \nu \rho \otimes \rho \otimes \sigma.
\end{equation}
According to \cite{Tad:israel}, Theorem 7.1, the representation 
$\rho \rtimes \delta(\nu \rho; \sigma)$ is reducible. 
Therefore, $\rho \rtimes \delta(\nu \rho; \sigma)=\tau + \tau'$,
where $\tau$ and $\tau'$ are irreducible tempered representations \cite{Gol}.
It follows from
$
  r_{M,G}(\rho \rtimes \delta(\nu \rho; \sigma) ) = 
   2 \, \rho \otimes \nu\rho \otimes \sigma
     + 2 \, \nu \rho \otimes \rho \otimes \sigma
$
that  
\begin{equation}\label{T}
 r_{M,G}(\tau) = \rho \otimes \nu\rho \otimes \sigma \quad \text{ and } \quad
 r_{M,G}(\tau') = \rho \otimes \nu\rho \otimes \sigma + 2 \, \nu \rho \otimes \rho \otimes \sigma.
\end{equation}
Further, $ r_{M,G}(\tau')$ exhausts appearance of 
$\nu \rho \otimes \rho \otimes \sigma$
in the Jacquet module of $\nu \rho \times \rho \rtimes \sigma$.
If we apply this information to (\ref{S}), we see that $\delta[\nu^{-1}\rho, \rho] \rtimes \sigma$
has two irreducible components, $\tau'$ and $\pi$. In addition,    
$r_{M,G}(\pi) = \rho \otimes \nu^{-1}\rho \otimes \sigma$.
We apply the Aubert involution and we obtain
$r_{M,G}(\hat{\pi}) = \rho \otimes \nu \rho \otimes \sigma$.
Equation (\ref{T}) tells us 
$
  \hat{\pi} = \tau \hookrightarrow \rho \rtimes \delta(\nu \rho; \sigma).
$
The $L$-parameter of $\hat{\pi}$ is
\[
    \phi \otimes S_3 \oplus  \phi \otimes S_1 \oplus
      \bigoplus_{i \in A} \phi_i \otimes S_1.
\]

Now, we will prove by induction on $m\geq 0$ that
\begin{equation}\label{U}
    \hat{\pi} = L_s(\delta[\nu^{-m}\rho, \nu^{-m+1}\rho], \dots ,
      \delta[\nu^{-2}\rho, \nu^{-1}\rho], \tau),
\end{equation}
where $\tau$ is, as above, the subrepresentation of $\rho \rtimes \delta(\nu \rho; \sigma)$
satisfying $r_{M,G}(\tau) = \rho \otimes \nu\rho \otimes \sigma$.
For $m=1$, (\ref{U}) follows from the first part of the proof.
Assume (\ref{U}) holds for $m$ and prove it holds for $m+1$.
Then $\pi = L_s(\delta[\nu^{-m-1}\rho, \nu^{m}\rho], \sigma)$.
We have to show
\begin{equation}\label{W}
    \hat{\pi} = L_s(\delta[\nu^{-m-1}\rho, \nu^{-m}\rho], \dots ,
      \delta[\nu^{-2}\rho, \nu^{-1}\rho], \tau).
\end{equation}
In the same way as in the proof of Theorem \ref{Th2}, we obtain
\[
   \pi  \hookrightarrow \nu^{m}\rho \times \nu^{m+1}\rho \times
    L_s(\delta[\nu^{-m}\rho, \nu^{m-1}\rho], \sigma).
\]
We apply the Aubert involution and the induction assumption to prove
that $\hat{\pi}$ is a component of
$
   \Pi= \nu^{-m}\rho \times \nu^{-m-1}\rho \times
  \delta[\nu^{-m}\rho, \nu^{-m+1}\rho] \times \cdots \times
      \delta[\nu^{-2}\rho, \nu^{-1}\rho] \rtimes \tau.
$
Let
\[
   \Pi_0 = \nu^{-m}\rho \otimes \nu^{-m-1}\rho \otimes
  \nu^{-m+1}\rho \otimes \nu^{-m}\rho \otimes \cdots \otimes
      \nu^{-1}\rho \otimes \nu^{-2}\rho \otimes \rho \otimes \nu\rho \otimes \sigma.
\]
It can be shown the multiplicity of $\Pi_0$ in $r_{M,G}(\Pi)$ is one.
Also, 
$\Pi_0 \leq r_{M,G}(\hat{\pi})$.
It follows that $\hat{\pi}$ is a subrepresentation of $\Pi$.
Since $\Pi$ is the representation induced from Langlands data, it has 
the unique subrepresentation. This proves (\ref{W}).
The $L$-parameter of (\ref{W}) is
\begin{equation}\label{V}
   \hat{\varphi} = (\bigoplus_{j = \frac{3}{2}}^{m+\frac{1}{2}} 
       |\cdot|^j  \phi \otimes S_2 \oplus |\cdot|^{-j}  \phi \otimes S_2)
     \oplus  \phi \otimes S_3 \oplus  \phi \otimes S_1 \oplus
      \bigoplus_{i \in A} \phi_i \otimes S_1.
\end{equation}
If we compare (\ref{V}) with (\ref{A}), we see that (\ref{V}) 
is not of the form  (\ref{A}). Indeed, the summand corresponding to
$j= \frac{1}{2}$ is missing in  (\ref{V}).
This completes the proof of the theorem.
\hfill  $\square$

\begin{corollary}\label{Cor}
Let $\rho$ be an irreducible unitary supercuspidal representation of $GL(n,F)$
and $\sigma$ an irreducible supercuspidal generic representation of $SO(2\ell+1,F)$.
Suppose  $\nu \rho \rtimes \sigma$ is reducible.
 \begin{enumerate}
   \item[{(i)}] If 
$
    \pi =L_s(\delta [\nu^{-m-\frac{1}{2}} \rho ,\nu^{m-\frac{1}{2}}\rho], \sigma),
$
then the Aubert involution of $\pi$ is equal to
\[
   \hat{\pi} 
= L_s(\nu^{-m-\frac{1}{2}} \rho, \nu^{-m+\frac{1}{2}} \rho,\nu^{-m+\frac{1}{2}} \rho, 
\dots, \nu^{-\frac{1}{2}} \rho,\nu^{-\frac{1}{2}} \rho, 
    \sigma). 
\]
\item[{(ii)}] If
$
    \pi =L_s(\delta [\nu^{-m} \rho ,\nu^{m-1}\rho], \sigma),
$
then
\[
   \hat{\pi} = L_s(\delta[\nu^{-m}\rho, \nu^{-m+1}\rho], \dots ,
      \delta[\nu^{-2}\rho, \nu^{-1}\rho], \tau),
\]
where $\tau$ is the subrepresentation of $\rho \rtimes \delta(\nu \rho; \sigma)$
satisfying $r_{M,G}(\tau) = \rho \otimes \nu\rho \otimes \sigma$.

\end{enumerate}

\end{corollary}

\noindent
{\bf Remark 6.1.} Let us consider the case $k=2$ of Theorem \ref{Th3}. 
Then
\[
    {\psi}= \phi \otimes S_{2} \otimes S_{2}
\oplus
     \bigoplus_{i \in A} \phi_i \otimes S_1 \otimes S_1.
\]
Let $\pi$ be the base point representation corresponding to $\phi_\psi$.
According to Corollary \ref{Cor}, $\hat{\pi} = \tau$, 
where $\tau$ is the subrepresentation of $\rho \rtimes \delta(\nu \rho; \sigma)$
satisfying $r_{M,G}(\tau) = \rho \otimes \nu\rho \otimes \sigma$.
The representation $\tau$ is tempered, but not generic. The generic representation
in the $L$-packet of $\tau$ is 
the subrepresentation $\tau'$ of $\rho \rtimes \delta(\nu \rho; \sigma)$
satisfying 
$
 r_{M,G}(\tau') = \rho \otimes \nu\rho \otimes \sigma + 2 \, \nu \rho \otimes \rho \otimes \sigma.
$
Therefore, in this example, the Aubert involution does not send a base point
 to a base point.
Set
\[
\psi' = \hat{\psi} = \phi \otimes S_3 \otimes S_1 \oplus  \phi \otimes S_1 \otimes S_1 \oplus
      \bigoplus_{i \in A} \phi_i \otimes S_1 \otimes S_1.
\]
It follows from \cite{BanZh} that
\[
\hat{\psi'} = \phi \otimes S_1 \otimes S_3 \oplus  \phi \otimes S_1 \otimes S_1 \oplus
      \bigoplus_{i \in A} \phi_i \otimes S_1 \otimes S_1,
\]
so $\psi'$ and $\hat{\psi'}$ are symmetric.

\end{document}